\documentclass[
a4paper,reqno,12pt]{amsart}

\usepackage{amsthm}
\usepackage{amsmath}
\usepackage{amssymb}

\begin{document}
\def\mmm{\setlength{\baselineskip}{3mm}}
\def\qq{\marginparsep plus \marginparwidth}
\def\pp{\qq plus \qq}
\def\VV{V_p(\M,a)}
\def\QQ{Q_p(\M,a)}
\def\RE{R(\varepsilon)}
\newcommand{\ox}{\overline{x}}
\newcommand{\mm}[1]{|#1|}
\newcommand{\myint}[3]{\int\limits_{#1}^{}\frac{#2}{#3}}
\newcommand{\xb}{x_a^{\bot}}
\newcommand{\xa}{x(m) - a}
\newcommand{\xt}{x_a^{\top}}
\newcommand{\xam}{x_a (m)}
\newcommand{\R}[1]{{\mathbb R}^{#1}}
\newcommand{\ma}[2]{M_{#1}(#2)}
\newcommand{\scal}[2]{{\mathcal  h} #1, #2{\mathcal  i}}
\newcommand{\DO}{\partial {\mathcal  O}_k(t)}
\newcommand{\MG}{{\mathcal  M}_g}
\newcommand{\M}{{\mathcal  M}}
\newcommand{\D}{{\mathcal  D}}
\newcommand{\N}{{\mathcal  N}}
\vspace*{4cm}
\begin{center}
\begin{bf}
FINITENESS OF THE NUMBER OF ENDS OF \\
MINIMAL SUBMANIFOLDS IN EUCLIDEAN SPACE\footnote{
Published in: \textit{Manuscr. Math.}, \textbf{82}(1994), no~1, 313-330}
\end{bf}
\end{center}
\vspace*{1 cm}
\begin{center}
Vladimir G. Tkachev\footnote{This paper was supported by Russian
Fundamental Research Foundation, project 93-011-176}
\end{center}

{
\begin{small}
\setlength{\baselineskip}{2mm}
We prove a version of the well-known Denjoy-Ahlfors theorem about the
number of asymptotic values of an entire function for properly immersed
minimal surfaces of arbitrary codimension in $\R{N}$.
The finiteness of the number of ends is proved for
minimal submanifolds with finite projective volume.
We show, as a corollary, that a minimal surface of codimension $n$ meeting
any $n$-plane passing through the origin in at most $k$ points has no more
$c(n,N) k$ ends.
\end{small}
}

\vspace*{0.5cm}
{\baselineskip = 5.5 mm
Let $x: M \rightarrow \R{n}$ be a proper minimal immersion of a
$p$-dimensional
orientiable manifold $M$. Then it is well-known that $M$ is necessarily
noncompact. The simplest topological invariant of such manifolds
is {\it the number of infinite points {\rm (}or ends
{\rm )}} of $M$, i.e. the smallest
integer $\ell (M)$ satisfying the following property: for every compact set
$F \subset M$ the number of the different components with noncompact closure
of $M \setminus F$ is less or equal to $\ell (M)$.

We say that a manifold $M$ (or the properly
immersed surface $\M = (M, x)$) is {\it manifold {\rm (}surface --
{\rm respectively}{\rm)} with finitely many ends} if
$\ell (M) < +\infty$.

These definitions agree with the usual ones for Riemannian surfaces
of finite type (see Example 1 below) and are related
to the noncompactness of the manifold.

In this paper we obtain some upper bounds for
$\ell (M)$ in terms of the  {\it projective volume} $V_p(\M )$
of $M$ and certain integral-geometric characteristics related to the
geometry in the large of minimal surfaces.

If ${\rm dim} M = 2$ and $M$ has finite total curvature $K(M)$,
R.Osser- man \cite{Os} (see also \cite{W}) proved that M is conformally
equivalent to a compact Riemann surface that has been punctured in a finite
number of points $\{m_1, m_2, \ldots, m_k\}$. In this case $\ell (M)$ is equal
to $k$. We observe, however, that the quantity $K(M)$ itself does not
represent any information about $\ell (M)$. Furthermore, there exist minimal
surfaces of finite topological type with $K(M) = -\infty$. For a detailled
discussion of these questions we refer to \cite{HM1}, \cite{HM2}.

The projective volume is one of the main tools in uniformization theory
and potential theory.
Using the special technique of estimating extremal lengths
in terms of a projective volume,
V.M.Miklyukov and the author in \cite{MT} showed
that a minimal surface $\M$ in $\R{3}$ has parabolic conformal type
provided that the generic number of points which  $\M$ has in common with
a line $L$ passing through a fixed point (possibly the infinitely far one) in
$\R{3}$ is uniformly bounded on $L$.
In particular, an upper bound for the projective volume
of the such surfaces was established.

In part 2 we prove that $\ell (M)$ is bounded by $c(p,n) V_p(\M )$.
We consider Theorem 2 as an extension of the Denjoy-Ahlfors theorem
about the number of asymptotic values (see the beautiful review
of A.Baernstein \cite{B}) to minimal submanifolds.
As a corollary we obtain in part 3
that a $p$-dimensional properly immersed minimal
surface meeting any $(n-p)$-plane, passing through origin,
in at most $k$ points  has no more than $c(p,n) k$ ends.
In particular, if a minimal hypersurface $\M$ is starlike with respect
to some point, then the number $\ell (M)$ is less than a constant
depending only on  ${\rm dim}\M$.

The results proved in this paper allow also to infer a parabolic
conformal type for properly immersed minimal submanifolds of
arbitrary codimension in the same way as in
\cite{MT}. We wish to mention also the paper \cite{WX} devoted to the
study of surfaces of hyperbolic type and \cite{Sch}
where results close to ours has been obtained.

\vspace*{0.6cm}
{\it I wish to thank  V.M.Miklyukov
for many useful discussions concerning the topic of this paper.
I also want to express particular appreciation to Professor Klaus
Steffen and the referee for many helpful suggestions
that greatly improved the presentation of this paper.}

\vspace*{1cm}

\noindent
{\bf 1. Some properties of the projective volume}}
\vspace*{0.7cm}

Let $a \in \R{n}$ and ${\mathcal  P}_a$ be the group consisting of
all conformal transformations preserving the set $\{ a, \infty
\}$, i.e ${\mathcal  P}_a$ is generated by the inversions: $x \to
\lambda~ (x~-~a~)~|~x~-~a~|^{-2}$ and the homotheties: $x \to
\lambda (x-a)$, where $\lambda $ is a positive factor.

Let $\M$ be a $p$-dimensional surface in $\R{n}$ and $B_a(R)$ be a
ball $\{x \in {\rm R}^n : \mm{x-a} < R \}$. We denote by  $M_a(R)$
the part of the surface ${\mathcal  M}$ inside $B_a (R)$ and
abbreviate $\xam = \xa$~, i.e. $\xam$ is a radius-vector of
$\mathcal M$ associated with  $a \in {\rm R}^n$. For given $a \in
\R{n}\setminus x(M)$ we define the following metric characteristic
of $\mathcal  M$ :

\begin{equation}
V_p(\M,a) = \limsup_{R \to \infty} \frac{1}{\ln R}\myint{M_a(R)}{1}
{\mm{\xam}^p}.
\label{defP}
\end{equation}

It is easy to see that $V_p(\M,a)$ is invariant under the action
of the group ${\mathcal  P}_a$. We call $V_p(\M,a)$ the {\it
projective volume} of $\M$.

Let $y^{\bot}(m)$ be the projection of $y$ on the normal space
to the surface $\M$ at a point $m$. Then we let
$$
Q_p(\M,a) = \myint{M}{\mm{\xb(m)}^2}{\mm{\xam}^{p+2}} =
\myint{M}{\mm{(\xa)^\bot}^2}{\mm{\xa}^{p+2}},
$$
and $ Q_p(\M,a) = +\infty$, if the last integral is divergent.
For $a \in x(M)$ we set
$$
Q_p(\M,a) = \lim_{\varepsilon \to 0} \myint{\mm{\xam} >
\varepsilon}{\mm{\xb (m)}^2}
{\mm{\xam}^{p+2}}.
$$

\bigskip
{\bf Theorem 1.} {\it Let $\mathcal  M$  be a  properly  immersed
$p$-dimensional minimal  surface in $\R{n}$ with compact boundary
$\Sigma$. Then the value $\VV$ does not depend on the choice of
$a\in\R{n}\setminus x(M)$. Moreover, the upper limit in {\rm
(\ref{defP})} can be replaced by a limit and
\begin{equation}
p\QQ = \VV + c(\Sigma;a),
\label{rr}
\end{equation}
where $c(\Sigma;a)$ is the finite constant such that $c(\varnothing ;a)=0$.
}

\bigskip

\rm
{\it Proof.} Let us assume $a\not\in x(M)$.
Denote $h= {\rm dist}(a,x(M))$, $r=\max_{m\in\Sigma}|x(m)-a|$
($r=0$, if $\Sigma = \varnothing$)
and $\varrho=\max\{ h;r\}$.
It is obvious from the properness of the immersion, that
 $\varrho > 0 $. Letting $f(m) = \mm{\xa}$ we have
$$
\nabla f(m) = \frac{\xt (m)}{\mm{\xam}}
$$
where  $(\; )^{\top}$ is the tangent part of the corresponding vector
and hence
$$
{\rm div} \frac{\xt (m)}{\mm{\xam}^p} =
\frac{1}{\mm{\xam}^p}\;{\rm div} (\xt (m))
- \frac{p}{\mm{\xam}^{p+1}}\;\scal{\nabla f}{\xt (m)}
$$
\begin{equation}
=\frac{p(\mm{\xam}^2 - \mm{\xt (m)}^2)}{\mm{\xam}^{p+2}} =
\frac{p\mm{\xb (m)}^2}{\mm{\xam}^{p+2}}.
\label{vova}
\end{equation}

Applying  Stokes'  formula  in  the   last   identity   over
$M_a(t,R)\equiv M_a(R)~\setminus~\overline{\ma{a}{t}} $
for $\varrho < t < R$  we have
$$
\frac{1}{R^p}\;\int\limits_{\partial \ma{a}{R}}{}\scal{\xt (m)}{\nu}
- \frac{1}{t^p}\;\int\limits_{\partial \ma{a}{t}}{}\scal{\xt (m)}{\nu} =
p\myint{M_a(t,R)}{\mm{\xb (m)}^2}{\mm{\xam}^{p+2}},
$$
where $\nu $ is the unit outward normal to the $t$-level set $\partial
 \ma{a}{t}$ of the
function $f(m)$. It is easy to see that for any regular
value  $t > \varrho$ of $f$ the normal $\nu$ is represented
on $\partial \ma{a}{t}$ by
$$
\nu (m) = \frac{\xt (m)}{\mm{\xt (m)}},
$$
and  thus  the last integral expression can be rewritten in the form
\begin{equation}
\frac{J(R)}{R^{p-1}} - \frac{J(t)}{t^{p-1}} =
p\myint{M_a(t,R)}{\mm{\xb (m)}^2}
{\mm{\xam}^{p+2}},
\label{jj}
\end{equation}
where $J(t) = t^{-1}\int_{\partial \ma{a}{t}}\mm{\xt (m)}$.
This relation is a well-known monotonicity formula for
minimal surfaces and (\ref{jj}) yields  the
increasing monotonicity of $J(t)t^{1-p}$.

On the other  hand  using  the  Kronrod-Federer  formula
(\cite{F}, Theorem 3.2.22 ), we obtain for  $R >R_1> \varrho $,
$$
\myint{M_a(R_1,R)}{1}{\mm{\xam}^p} =
\myint{M_a(R_1,R)}{\mm{\xt (m)}^2}{\mm{\xam}^{p+2}} +
\myint{M_a(R_1,R)}{\mm{\xb (m)}^2}{\mm{\xam}^{p+2}}
$$
$$
= \int\limits_{R_1}^{R}\frac{dt}{t^{p+2}}
\int\limits_{\partial M_a(t)}\frac{\mm{\xt (m)}^2}
{\mm{\nabla f}} +\myint{M_a(R_1,R)}{\mm{\xb (m)}^2}{\mm{\xam}^{p+2}}
$$
$$
= \int\limits_{R_1}^{R}\frac{J(t)}{t^{p-1}}\frac{dt}{t} +
\myint{M_a(R_1,R)}{\mm{\xb (m)}^2}{\mm{\xam}^{p+2}}
$$
\begin{equation}
= \int\limits_{R_1}^{R}\frac{J(t)}{t^{p-1}}\frac{dt}{t}\, + \,
\frac{1}{p}\left(
\frac{J(R)}{R^{p-1}} - \frac{J(R_1)}{R_1^{p-1}}
\right).
\label{star}
\end{equation}
Now we notice that the increasing monotonicity of $J(t)t^{1-p}$
yields immediately
\begin{equation}
\lim_{R\to+\infty}\frac{J(R)}{R^{p-1}}= \VV,
\label{ss}
\end{equation}
and consequently,
\begin{equation}
\frac{J(R)}{R^{p-1}}\leq \VV ,\quad \mbox{for} \quad R>\varrho.
\label{ee}
\end{equation}
Then we deduce from (\ref{ss}) and (\ref{star}) that the upper
limit in (\ref{defP}) can be replaced on a limit.

Now we show that really $V_p(\M,a) $ does not depend on $a$.
Let $b$ be a point in $\R{n}\setminus x(M)$ and $\delta = |b-a|$.
Then for every $0<\varepsilon <1$ there exist $R(\varepsilon )>\delta$
and $c(\varepsilon)<\infty$ (both independent of $R$) such that for
$R>R(\varepsilon)$ one has
$$
\left |\myint{M_b(R)}{1}{\mm{x-b}^p} - \myint{M_b(R)}{1}{\mm{x-a}^p}\right |
 \leq \int\limits_{M_b(R(\varepsilon))}
 \left|\frac{1}{|x-a|^p} -\frac{1}{|x-b|^p}\right|
$$
\begin{equation}
+p\delta\int\limits_{M_b(R(\varepsilon),R)}\frac{\mm{x-a}^{p-1} +
\mm{x-b}^{p-1}}{\mm{x-b}^p\mm{x-a}^p}
\leq
c(\varepsilon)+\varepsilon \myint{M_b(R)}{1}{\mm{x-b}^{p}},
\label{o1}
\end{equation}
(choosing $\RE>>1$ such that
$\frac{p\delta}{\RE -\delta}\left [
\left( \frac{\RE}{\RE -\delta}\right )^{p-1}+1\right ]\leq \varepsilon$).

Moreover, from the obvious inclusions
$$
M_a(R-\delta) \subset M_b(R) \subset M_a(R+\delta),
$$
for $ R > \delta$, we have
$$
\myint{M_a(R-\delta)}{1}{\mm{x-b}^p} \leq \myint{M_b(R)}{1}{\mm{x-b}^p}
\leq \myint{M_a(R+\delta)}{1}{\mm{x-b}^p}.
$$
From this and (\ref{o1}),
it then follows that
\begin{eqnarray}
\frac{1}{1+\varepsilon}\myint{M_a(R-\delta)}{1}{|x-a|^p}-
\frac{c(\varepsilon)}{1+\varepsilon}
&
\leq
&\myint{M_b(R)}{1}{|x-b|^p}
\nonumber\\
&\leq&\frac{1}{1-\varepsilon}\myint{M_a(R+\delta)}{1}{|x-a|^p}+
\frac{c(\varepsilon)}{1-\varepsilon}
\nonumber
\end{eqnarray}
and dividing by $\ln R$ one infers
$$
\frac{1}{1+\varepsilon}V_p(\M,a)\leq
V_p(\M,b)\leq
\frac{1}{1-\varepsilon}V_p(\M,a)
$$
letting $R\to\infty$. In view of arbitrariness of $\varepsilon$
this implies $V_p(\M,a)=V_p(\M,b)$ whether both quantities are
finite or not.

Now integrating (\ref{vova}) over $M_a(R)$, we obtain for $R>\varrho$
$$
\frac{1}{R^p}\;\int\limits_{\partial \ma{a}{R}}{}\scal{\xt (m)}{\nu}
- \int\limits_{\Sigma}{}\frac{\scal{\xt (m)}{\nu}}{|x_a(m)|^p} =
p\myint{M_a(R)}{\mm{\xb (m)}^2}{\mm{\xam}^{p+2}}.
$$
Let us denote by $c(\Sigma;a)$ the second integral in the above equation.
Then we have
$$
\frac{J(R)}{R^{p-1}}-c(\Sigma;a)=
p\myint{M_a(R)}{\mm{\xb (m)}^2}{\mm{\xam}^{p+2}}.
$$
Leting $R\to\infty$ and using the equality (\ref{ss}), we complete
the proof of Theorem 1.

\smallskip
We denote by $ a\#\mathcal  M$ the multiplicity of the immersion
$x : M \rightarrow \R{n}$ at the point $ a~\in~\R{n} $~, i.e. the
cardinal number of the preimage $x^{-1}(a\cap x(M))$.

The next property of $V_p(\M,a) $ and   $Q_p(\M,a)$ shows that
these quantities are {\it conformal invariants} of minimal
surfaces  in the sense that $V_p(g\circ \M; a) = V_p(\M ; 0)$ for
all $g \in {\mathcal  P}_a$, and $ a \in \R{n}\setminus x(M)$.

\bigskip
{\bf Corollary 1.} {\it Let $\M\subset \R{n}$  be a  properly
immersed  minimal  surface without boundary, $\dim {\mathcal  M} =
p$. Then both values $Q_p(\M,a)$ and $V_p(\M,a)$ do not depend on
choice of $a \in \R{n}\setminus x(M)$ and for all $a\in \R{n}$ we
have
$$
\frac{1}{p}\; V_p(\M,a) = Q_p(\M,a) - \omega_p (a\# {\mathcal
M}),
$$
where $\omega_p$  is the $(p - 1)$-dimensional Hausdorff measure of
unit sphere} $S^{p-1}$.
\rm

\bigskip
{\it Proof.} We observe that the first assertion of Corollary 1  follows
immediately from $\Sigma =\varnothing$ and (\ref{rr}).

Let us now consider $a \in x(M)$ so that $q = a\#\mathcal  M$ is a
positive integer. Reasoning similarly as above we get that
$J(t)t^{1-p}$ is a  positive monotonic function for $t \to +0$.
Consequently, there exists
$$
\mu = \lim_{t \to +0}\frac{J(t)}{t^{p-1}}.
$$
We consider any preimage $m_k \in x^{-1}(a)$. Let ${\mathcal
O}_k(t)$ be  an  open component of $M_a(t)$ which contains $m_k$.
It is clear  that  for sufficiently small $t > 0$ the sets
${\mathcal O}_k(t)$ are nonintersecting for  all $k \leq q$. Then
from (\ref{jj}) we have
\begin{equation}
\frac{J(R)}{R^{p-1}} - \mu =
\; p\myint{M_a(R)}{\mm{x_a(m)^{\bot}}^2}{\mm{x_a(m)}^{p+2}}.
\label{jmu}
\end{equation}
But in virtue of the regularity of the immersion $x(m)$,
$$
\lim_{t\to +0}\; \frac{1}{t} \sup_{m\in {\mathcal  O}_k(t)}\mm{\xb
(m)} \quad = 0.
$$

It follows that
\begin{equation}
\lim_{t\to +0} \frac{1}{t^{p-1}} \myint{\partial {\mathcal
O}_k(t)}{|\xt (m)|} {\mm{x_a(m)}} = \lim_{t\to +0} \frac{{\rm
meas}_{p-1}(\partial{\mathcal  O}_k(t))}{t^{p-1}} = \omega_p,
\label{limmes}
\end{equation}
and taking into consideration that for small $t > 0$
$$
\partial M_a(t) = \bigcup_{k=1}^{q} \DO,
$$
we obtain $\mu = q\;\omega_p$.

Repeating the above arguments we conclude
$$
\lim_{R \to\infty}\frac{1}{\ln R}\myint{M_a(R)}{1}{\mm{\xam}^p} =
\lim_{R\to +\infty}\frac{J(R)}{R^{p-1}} \quad,
$$
and by (\ref{jmu}) and (\ref{limmes}) Corollary 1 is proved.

From  now on we  write $Q_p(\M )$ and $V_p(\M )$ instead of
$Q(\M,a)$ and $V_p(\M,a)$ respectively, if $a\not\in x({\mathcal
M})$.

{\it Example} 1. Let $\MG $ be a compact orientiable Riemannian
surface $\MG $ of genus $g \geq 0$, and $m_1, m_2, \ldots m_l \in
\MG$. Let $\zeta$ be a holomorphic $1$-form on $\MG$ and $h: \MG
\rightarrow {\mathbb C}\cup\{ \infty \}$ a meromorphic function.
Then due to \cite{Os} the vector valued $1$-form
$$
\Phi = (\Phi_1, \Phi_2, \Phi_3 )^t =
\left ((1 - h^2)\zeta ;\, i(1 + h^2)\zeta ;\, 2h\zeta\right )^t
$$
gives a conformal minimal immersion
$$
X(m) = {\rm Re} \int\limits_{m_0}^{m} \Phi
$$
which is well-defined on ${\mathcal  M}^\star_g =
\MG\setminus\{m_1, m_2,\ldots,m_l\}$ and regular, provided
\begin{enumerate}
\item
\it
No component of $\Phi $  has a real period on $\MG$;
\rm
\item
\it
The poles $\{m_1, m_2,\ldots,m_l\}$  of $h$ coincide with zeros of $\zeta$ and
the order of a pole $m_k$ of $h$ is precisely the order of the corresponding
zero of $\zeta$.
\end{enumerate}
\rm
It is well-known in the case of finite total curvature
that the asymptotic behaviour of $X(m)$  in the neighbourhood of $m_k$ is
either of flat or catenoid type \cite{HM2}.
In both cases the quantity $Q_p(\M )$ and, consequently, $V_p(\M )$
can be calculated directly and we have
$$
V_p(\M ) \equiv 2 Q_p(\M )= 2 \pi l \,.
$$
We observe that \it the characteristic $Q(\M )$ does not depend on
the genus $g$ of $\MG$, and describes only the noncompactness
nature of ${\mathcal  M}^\star_g$.

\rm
{\it Remark} 1. It would be interesting to know in analogy
with the case of finite Gaussian curvature above, {\it whether the
set of possible values of the quantity $Q(\M )$ is discrete}. It follows
from the above example that this is true for two-dimensional minimal
surfaces of finite topology.

\vspace*{1cm}
\noindent
{ \bf 2. The estimate
for the number of ends of minimal submanifolds}
\vspace*{0.7cm}

In this section we give a geometric application of the above invariants.

\bigskip
{\bf Theorem 2.}
\it
 Let $\mathcal  M$ be  a  properly  immersed  $p$-dimensional
minimal surface in $\R{n}$ with compact boundary $\Sigma$, having
finite projective volume $V_p(\M )$. Then $\mathcal  M$ is a
surface with finitely many ends and
$$
\ell (M) \leq \frac{2^p}{\omega_p} V_p(\M ).
$$

\rm
\bigskip
The proof of the theorem is based on  the  next  auxiliary  assertion.

\bigskip
{\bf Lemma 1.} \it Let $\mathcal  D$ be a connected
$p$-dimensional minimal surface with boundary $\partial{\mathcal
D} \subset \partial B_0(R_1)\cup\partial B_0(R_2), \quad R_2 > R_1
> 0.$ Then
\begin{equation}
{\rm meas}_p {\mathcal  D} \geq \frac{\omega_p}{p}\,\left
(\frac{R_2-R_1}{2}\right )^p. \label{lemD}
\end{equation}
\rm

\bigskip
{\it Proof} of the lemma. We consider first the case when
$\mathcal  D$  is  a compact minimal submanifold such that $0 \in
\mathcal  D$ and $\partial{\mathcal  D} \subset \partial B_0(R) $.
Let
$$
A(t) = {\rm meas}_p({\mathcal  D}\cap B_0(t))\,.
$$
Then using
$$
{\rm div}x^{\top} (m) = \sum_{i=1}^{n} {\rm div}x_ie_i^\top =
\sum_{i=1}^{n} \mm{e_i^\top}^2 = p
$$
we shall have after integration
\begin{equation}
p A(t) = \int\limits_{{\mathcal  D}\cap B_0(t)} {\rm div}x^{\top}
(m) = t \myint{{\mathcal  D}\cap \partial B_0(t)}{\mm{x^\top
(m)}}{\mm{x(m)}} =t\, J(t), \label{At}
\end{equation}
that the function
$$
\frac{p\, A(t)}{t^p} = \frac{J(t)}{t^{p-1}}
$$
is an increasing one. Moreover,
$$
\lim_{t\to +0} \frac{p\, A(t)}{t^p} = \lim_{t\to
+0}\frac{J(t)}{t^{p-1}} = \omega_p\, \cdot (0\#{\mathcal  D}),
$$
and consequently, for all $t > 0$ we have
\begin{equation}
\frac{{\rm meas}_p({\mathcal  D}\cap B_0(t))}{t^p \omega_p} \geq
\lim_{t\to +0} \frac{A(t)}{t^p} = \frac{1}{p} \,(0\#{\mathcal  D})
\label{mesA}
\end{equation}
and the first case of the lemma is proved.

Let  us  now  assume  that  $\partial {\mathcal  D} \subset
\partial B_0(R_1) \cup\partial B_0(R_2)$. Put $R = \frac{1}{2}(R_1
+ R_2)$. We observe that the set $ {\mathcal  D}\cap\partial
B_0(R)$ is not empty by virtue of the connectivity $\mathcal  D$,
and we let $a$ be  any  point  in ${\mathcal  D}\cap\partial
B_0(R)$. Then for $\D_1 = \D \cap B_0(r)$ we have
$$
\partial\D_1 \subset \partial B_0(r)
$$
for $r = \frac{1}{2} (R_2-R_1)$.
In view of (\ref{mesA}) and the inclusion $\D_1 \subset \D$,
the  above  inclusion  implies
(\ref{lemD}) and thus the proof the lemma is concluded.

\medskip
{\it Remark} 2. We note that Lemma 1 can be also obtained from the
general result of W.K.Allard \cite{A}.

\medskip
{\it Proof} of Theorem 2.
Without  loss  of generality we can arrange that $0 \not\in
\overline{\mathcal  M}$ and, by Theorem 1, $V_p(\M )= V_p(\M,0)$.
We fix a sufficiently large regular value $R > 0$ of $f(m)=|x(m)|$
such that $\Sigma \subset B_0(R)$.

Let $\D_1, \ldots, \D_k \ldots$ be the open components of
$M~\setminus~\overline{M_0(R)}$.
Notice that $x(\partial \D_k)\subset \partial B_0(R)$
and $\Delta~f(m)~\geq~0$.
Then the maximum principle implies that the $\D_k$ are domains with noncompact
closure.
Moreover, it follows from the regularity of $R$ that
the number $l=l(R)$ of components $\D_k$ is finite,
and it is nondecreasing with respect to $R$.
Put for $t > R$,
$$
J_k(t) = \frac{1}{t} \int\limits_{\partial B_0(t)\cap\D_k}\mm{x^\top (m)}.
$$
Then reasoning  similarly as in the  proof  of Theorem 1 we arrive at the
inequality
\begin{equation}
\sum_{k=1}^{l} J_k(t) = J(t) \;\leq \;
V_p(\M )t^{p-1}.
\label{sumj}
\end{equation}
On the other hand, applying Lemma 1, we have
$$
{\rm meas}_p\{m\in\D_k :\; |x(m)|<t\} \geq \frac{\omega_p (t-R)^p}{p\, 2^p},
$$
and after summing over all $k\leq l$  we obtain that
$$
{\rm meas}_p(M_0(t)\setminus M_0(R)) \geq
\frac{l \omega_p}{p 2^p} (t-R)^p.
$$
Using (\ref{At}) and (\ref{sumj}) we have the sequence of inequalities:
$$
\frac{l \omega_p}{p 2^p} (t-R)^p \leq {\rm meas}_p(M_0(t)\setminus M_0(R))
\leq {\rm meas}_p(M_0(t)) =
$$
$$
= \frac{tJ(t)}{p} \leq\frac{V_p(\M )t^p}{p},
$$
and after dividing by $t^p$ and letting $t \to \infty$, we obtain
$$
l(R)=l \leq \quad \frac{V_p(\M )2^p}{\omega_p}.
$$
Next,  from  the  fact  that  the  integer-valued  function  $l(t)$
is nondecreasing, we conclude that it  is  stabilized,  i.e.
$l(R) \equiv \rm const$ for sufficient large $R$.

Let $F\subset M$ be an arbitrary compact subset. Using again the
the maximum principle and the properness of immersion
that the number of components with noncompact closure
of $M\setminus F$ with noncompact closure is  a nondecreasing
function of the compact set $F$. Therefore,
$\ell (M)\equiv \lim_{t\to\infty}l(t)$ and the theorem is proved.

\medskip
{\bf Corollary 2.}
\it
Let $\M$ be  a  properly  immersed  $p$-dimensional
minimal surface without boundary  having finite  projective
volume. Then $\M$ is a surface with finitely many ends and
$$
\ell(M) \quad \leq \quad  \frac{Q_p(\M ) 2^pp}{\omega_p}.
$$
\rm


\bigskip
\noindent
{\bf 3. The bounded integral-geometric averages and the finite-}

\noindent
{\bf ness of the number of ends of minimal submanifolds}

\bigskip

In this section we discuss certain sufficient conditions
for the finiteness of  the  projective  volume  for  minimal
submanifolds with arbitrary codimension.

Suppose first  that $\M$  is  a  hypersurface  in $\R{n}$.
Then
specifying a point $b \in \R{n}\setminus x(M)$ we can introduce  {\it the
counting  function}
$\N (e, b)$ for the multiplicity  of the radial projection relative to $b$,
setting for any unit direction $e \in\R{n}$
$$
\N (e, b) = \sum_{a\in L_b(e)}a\#\M \;  \equiv \;\# x^{-1}(L_b(e)\cap x(M)),
$$
where $L_b(e)$ is a ray with the origin  at $b$   directed  as $e$.  The
number $\N (e, b)$ can be interpreted  as  the  multiplicity  of  the
covering
\begin{equation}
\pi_b : \;\M\rightarrow S^{n-1},\quad
\pi_b(y) = \frac{y-b}{\mm{y-b}},
\label{pib}
\end{equation}
at a point $e$.

If ${\rm codim}\M > 1$, then the image of $\M$  after  projection  (\ref{pib})
is a null-measure subset in $S^{n-1}$ and the second definition  of
$\N (e, b)$ is meaningless. Therefore we give the following generalization
of  the first definition.
\par
Let $G_n^p(b) $ be the Grassman manifold of all nonoriented
$(n-p)$-dimensional planes $\gamma $  passing through $b$. Then $G_n^p(b)$
can be equipped with the unique  Haar measure $d\gamma $  which is invariant
under the  action of the motion subgroup preserving $b$, and
normalized by
$$
\int\limits_{G_n^p(b)} d\gamma \; = \; 1.
$$
Let $R > 0$. By  Sard's  theorem we know  that  for
$d\gamma $-almost all planes $\gamma \in G_n^p(b)$  the set of the preimages
$x^{-1}(x(M) \cap\gamma\cap B_b(R))$ is  a discrete one. Put
$$
\N (b,\gamma; R) = \; \# x^{-1}(x(M) \cap\gamma\cap B_b(R)),
$$
- the cardinality of the corresponding set.
The quantity
$$
\N (b; R) = \int\limits_{G_n^p(b)} \N (b,\gamma; R)d\gamma
$$
can now be interpreted as "the average multiplicity"  of
the intersection of $(n-p)$-dimensional planes with  the  part  of $\M$
distant from $b$ not further than $R$. Moreover, $\N (b; R)$ is an
increasing function of $R$ and hence there exists a finite or infinite limit
$$
\N (b) = \; \lim_{R \to\infty}\N (b; R).
$$

\medskip
\it
{\bf Lemma 2.} Let $\M$  be  a  $p$-dimensional  properly  immersed
minimal surface in $\R{n}$ without boundary  and $b \not\in \M$. Then
\begin{equation}
Q(\M ) \leq \frac{1}{2}\, \N (b)\omega_{p+1},
\label{qun}
\end{equation}
where $\omega_{p+1}$  is the $p$-dimensional Hausdorff measure of
unit sphere $S^p$.

\rm\medskip\par
{\it Proof.} Without loss of generality we can assume
that $b = 0$~. We specify $R > 0$ and denote as above
$$
M_0(R) = \{ m \in M : \mm{x(m)} < R \}.
$$
We consider the composition
$$
\sigma : \M \stackrel{x}{\rightarrow} \R{n}\setminus\{ 0 \}
\stackrel{\pi}{\rightarrow} S^{n-1},
$$
where $\pi $ is defined as in (\ref{pib}) with $b = 0$. In order  to  find
the
Jacobian ${\rm det} (d\sigma ) $ of the map $\sigma $  at $m$  we observe that
$$
d\sigma_m = d\pi_{x(m)} \circ dx_m : T_mM \rightarrow T_{\sigma (m)}S^{n-1}.
$$
By direct calculation one can show that
$$
d\pi_a(X) = \frac{X - \pi (a)\scal{X}{\pi (a)}}{\mm{a}},
$$
for all $a \in\; \R{n}\setminus\{0\}$ and $X \in T_a\R{n}$. Hence for any
$Y \in T_mM$
$$
d\sigma_m(Y) = \frac{Y - \ox (m)\scal{Y}{\ox (m)}}{\mm{x(m)}}
$$
where $\ox (m) = x(m)/\mm{x(m)} $ and we identify $Y$ with $dx_m(Y)$
and $T_mM$ with a subspace of  $T_{x(m)}\R{n}\cong\R{n}$
through the isometry $dx_m$.
Choose an orthonormal basis $Y_1,\ldots,Y_p$ in $T_mM$.
We then have
$$
{\rm det}^2 (d\sigma_m) = \scal{w}{w}
$$
 where
$$
w = d\sigma_m(Y_1)\wedge d\sigma_m(Y_2)\wedge\ldots\wedge d\sigma_m(Y_p) =
$$
$$
= \mm{x}^{-p}\,\bigl ( Y_1 - \ox \scal{Y_1}{\ox }\bigr )\wedge\ldots\wedge
\bigl ( Y_p - \ox \scal{Y_p}{\ox }\bigr )
$$
$$
= \mm{x}^{-p} ( Y_1\wedge Y_2\wedge\ldots\wedge Y_p -
 \sum_{i=1}^{p}Y_1\wedge\ldots\wedge Y_{i-1}\wedge\ox \wedge Y_{i+1}\wedge
\ldots\wedge Y_p\scal{\ox }{Y_i} )
$$
$$
=\frac{ Y_1\wedge Y_2\wedge\ldots\wedge Y_p}{\mm{x}^p}
\left( 1 - \sum_{i=1}^{p} \scal{Y_i}{\ox (m)}^2\right)
$$
$$
 - \sum_{i=1}^{p}Y_1\wedge\ldots\wedge Y_{i-1}\wedge\ox^{\bot}
 \wedge Y_{i+1}\wedge \ldots\wedge Y_p\scal{\ox }{Y_i} ,
$$
and consequently,
$$
{\rm det}^2 (d\sigma_m) =
\frac{\mm{\ox^\bot (m)}^2}{\mm{x(m)}^{2p}}.
$$
Thus we obtain
\begin{equation}
|{\rm det} (d\sigma_m)| =  \frac{\mm{\ox^\bot (m)}}{\mm{x(m)}^{p}} =
 \frac{\mm{x^\bot (m)}}{\mm{x(m)}^{p+1}},
\label{det}
\end{equation}
the required expression for the Jacobian  of $d\sigma_m$.  By  the  change
of coordinates formula we obtain from (\ref{det})
$$
\myint{M(R)}{\mm{x^\bot (m)}^2}{\mm{x(m)}^{p+2}} \leq
\myint{M(R)}{\mm{x^\bot (m)}}{\mm{x(m)}^{p+1}}
$$
\begin{equation}
= \int\limits_{M(R)}\left |{\rm det} (d\sigma_m)\right | =
\int\limits_{\sigma(M(R))} \chi (s) d{\mathcal  H}^p(s),
\label{inthi}
\end{equation}
where $\chi (s)$ is the cardinality of the preimage $\sigma
^{-1}(s) \cap M(R)$ for the given  $s \in \sigma (M(R)) \subset
S^{n-1}$ and ${\mathcal  H}^p$ is the  corresponding Hausdorff
measure on $\sigma (M(R))$. According to the theorem of  Federer
(\cite{F}, Theorem 3.2.48), we conclude that for every ${\mathcal
H}^p$-measurable and  $({\mathcal  H}^p,p)$-rectifiable set  $F
\subset S^{n-1}$ and positive summable function $f$ on $F$
\begin{equation}
\int\limits_{s\in F} f(s) d{\mathcal  H}^p(s) =
\frac{\omega_{p+1}}{2}\int\limits_{\gamma\in G_n^p(0)}
f^{\ast}(\gamma\cap F) d\gamma, \label{federer}
\end{equation}
where $ f^{\ast}(\gamma\cap F) = \sum_{s\in \gamma\cap F} f(s) $
is  well-defined  function  for
$d\gamma $-almost all planes $\gamma \in G_n^p(0) $. Then it follows from
(\ref{inthi}) and (\ref{federer})
$$
\myint{M(R)}{\mm{x^\bot (m)}^2}{\mm{x(m)}^{p+2}} \leq
\frac{\omega_{p+1}}{2}\int\limits_{\gamma\in G_n^p(0)}
 \#\sigma^{-1}[\gamma\cap \sigma (M(R))] d\gamma =
\frac{\omega_{p+1}}{2} \N (0, R),
$$
and taking $R \to \infty$ we arrive at the required estimate (\ref{qun}).
\par
Thus using the previous lemma and Corollary 2, we have

\medskip\par\it
{\bf Corollary 3.} Let $\M$ be  a  properly  immersed  $p$-dimensional
minimal surface in $\R{n}$ without boundary.
Suppose that for some point $b \in \R{n}$ the
cardinality of the set of intersection points (taking into account
multiplicity) of any $\gamma \in G_n^p(b)$  and $x(M)$ does not exceed
$k$. Then $M$  is a manifold with finitely many ends and
$$
\ell(M) \leq  k c_p ,
$$
where
$$
c_p = 2^{p-1}(p+1)\sqrt{\pi}\Gamma (\frac{p+2}{2})\Gamma^{-1}(\frac{p+3}{2})=
\frac{2^{p-1} p\omega_{p+1}}{\omega_{p}}
$$
and $\Gamma$ is the  Euler gamma-function and $\omega_{p+1}$
is as in Lemma 2.
\rm

\bigskip\par\it
{\bf Corollary 4.} Let $\M$ be a properly embedded $p$-dimensional
minimal hypersurface without boundary. Assume that $\M$ is starlike with respect
to some point in $\R{p+1}$. Then  the number of ends $\ell(M)$ satisfies
$$
\ell(M) \leq 2 c_p,
$$
where the constant $c_p$ is from the previous lemma.

{\begin{small}
\baselineskip=2mm
\rm

\end{small}
}

\end{document}